\documentclass[12pt]{amsart}
\usepackage{amsmath,amssymb,amsfonts,amsthm,amsopn}
\usepackage{graphics}

  \def\Spnr{Sp(d,\R)}
  \def\Gltwonr{GL(2d,\R)}
  \newcommand\mc[1]{{\mathcal{#1}}}

\newcommand{\modsp}{modulation space}

\newtheorem{tm}{Theorem}[section]
\newtheorem{lemma}[tm]{Lemma}

\newtheorem{theorem}{Theorem}[section]
\newtheorem{corollary}[theorem]{Corollary}

\newtheorem{proposition}[theorem]{Proposition}
\newtheorem{remark}[theorem]{Remark}

\newcommand{\beqa}{\begin{eqnarray*}}
\newcommand{\eeqa}{\end{eqnarray*}}

\newcommand{\field}[1]{\mathbb{#1}}
\newcommand{\bR}{\field{R}}        
\def\la{\lambda}

 \def\cF{\mathcal{F}}              
 \def\cS{\mathcal{S}}
 \def\cD{\mathcal{D}}

 \def\cU{\mathcal{U}}

 \def\cC{\mathcal{C}}

\def\a{\aleph}

\def\rd{\bR^d}

\def\rdd{{\bR^{2d}}}

\def\lrd{L^2(\rd)}

\def\intrd{\int_{\rd}}
\def\intrdd{\int_{\rdd}}

\def\R{\right)}

\def\<{\left<}
\def\>{\right>}

\def\inv{^{-1}}

\def\mv1{M_v^1}

\def\mpq{M^{p,q}}

\def\phas{(x,\xi )}

\def\o{\xi}
\def\a{\alpha}
\def\b{\beta}

\def\ZZ{\mathbb{Z}}

\def\R{\mathbb{R}}
\def\Ren{\mathbb{R}^d}
\def\Renn{\mathbb{R}^{2d}}

\def\sch{\mathcal{S}}

\def\Fur{\mathcal{F}}

\def\f{\varphi}

\def\Sn2{S_{2}(L^{2}(\Ren))}
\def\S1{S_{1}(L^{2}(\Ren))}
\def\sig00{\sigma_{0,0}}

\def\la{\langle}
\def\ra{\rangle}

\def\spdr{{\mathfrak {sp}}(d,\R)}

\newcommand{\A}{\mathcal{A}}

\newcommand{\fpq}{W(\mathcal{F}L^p,L^q)}
\newcommand{\fqp}{W(\mathcal{F}L^q,L^p)}
\newcommand{\fui}{W(\mathcal{F}L^1,L^\infty)}
\newcommand{\fiu}{W(\mathcal{F}L^\infty,L^1)}
\newcommand{\frrp}{W(\mathcal{F}L^r,L^{r'})}
\newcommand{\frpr}{W(\mathcal{F}L^{r'},L^r)}
\newcommand{\mil}{\min\{1,|\lambda_j|\}}
\newcommand{\mal}{\max\{1,|\lambda_j|\}}

\begin{document}

\title[Metaplectic representation
and applications]{Metaplectic representation on Wiener amalgam
spaces and applications to the Schr\"odinger equation}

\author{Elena Cordero and Fabio Nicola}
\address{Department of Mathematics,  University of Torino, Italy}
\address{Dipartimento di Matematica, Politecnico di Torino, Italy}
\email{elena.cordero@unito.it}
\email{fabio.nicola@polito.it}

\date{}

\begin{abstract}
We study the action of metaplectic operators on Wiener amalgam
spaces, giving upper bounds for their norms. As an application, we
obtain new fixed-time estimates in these spaces for Schr\"odinger
equations with general quadratic Hamiltonians and Strichartz
estimates for the Schr\"odinger equation with potentials
$V(x)=\pm|x|^2$.
\end{abstract}

\maketitle

\section{Introduction}
The  Wiener amalgam spaces
 were introduced by
 Feichtinger \cite{feichtinger80} in
 1980 and soon they revealed
 to be, together with the closely related modulation spaces,
 the natural framework
 for the Time-Frequency
 Analysis; see e.g. \cite{feichtinger83,
feichtinger90,fournier-stewart85,Fei98}
and Gr\"ochenig's book
 \cite{grochenig}. These
 spaces are modeled on the
 $L^p$ spaces but they turn
 out to be much more flexible,
 since they control the local
 regularity of a function and its decay at
 infinity separately. For
 example, the Wiener amalgam space
 $W(B,L^q)$, $B=L^p$ or $B=\Fur
 L^p,$ etc.,  consists of
 functions which locally have
 the regularity of a function
 in $B$ but globally
  display a $L^q$ decay. \par
  In this paper we focus our attention on the action of
   the metaplectic
  representation on Wiener amalgam
  spaces. The metaplectic representation
  $\mu:Sp(d,\R)\to \cU(L^2(\R^d))$ of the
  symplectic group $Sp(d,\R)$
  (see the subsequent  Section 2 and \cite{folland} for details), was first constructed by Segal and
  Shale \cite{Segal, Shale} in the framework of quantum mechanics
  (though on the  algebra level the first construction is due to van
  Hove \cite{Hove}) and by Weil \cite{Weil} in number theory. Since
  then, the metaplectic representation has attracted the attention of
  many people in different fields of mathematics and physics.
  In particular, we highlight the applications in the framework of
  reproducing formulae and wavelet theory \cite{aefk1},  frame theory
  \cite{norbert}, quantum mechanics \cite{degosson}
      and  PDE's \cite{Hormander3,Hor95}.

Fix a test function $g\in\cC_0^\infty$ and $1\leq p,q\leq\infty$.
Then, the {\it Wiener amalgam space} $W(\Fur L^p,L^q)$ with local
component $\Fur L^p$ and global component  $L^q$ is defined as the
space of all functions/tempered distributions such that
$$\|f\|_{W(\Fur L^p,L^q)}: = \|\| fT_x g\|_{\Fur
L^p}\|_{L^q_x}<\infty,
$$
where $T_x g(t):=g(t-x)$. To give a flavor of the type of results:
\vskip0.3truecm {\it If $1\leq p\leq
  q\leq\infty$ and $\A=\begin{pmatrix} A&B\\
  C&D\end{pmatrix}\in Sp(d,\R)$, with
  $\det B\not=0$, then  the metaplectic operator $\mu(\A)$ is a
  continuous mapping from $W(\Fur L^q,L^p)$ into
  $W(\Fur L^p,L^q)$, that is
$$
\|\mu(\mathcal{A})f\|_{\fpq}\leq
\a(\mathcal{A},p,q)\|f\|_{\fqp}.
$$}
 The norm upper bound $\a=\a(\mathcal{A},p,q)$
 is explicitly expressed in terms of the matrix
 $\A$ and the indices $p,q$ (see Theorems \ref{t1} and \ref{ttt1}).

  This analysis generalize the basic result \cite{feichtinger83}:
\vskip0.3truecm
  \emph{The Fourier
 transform $\cF$ is a continuous mapping between  $W(\Fur L^q,L^p)$ and
  $W(\Fur L^p,L^q)$
   if (and only
  if) $1\leq p\leq
  q\leq\infty$ }.
\vskip0.3truecm

  Indeed, the
  Fourier transform $\cF$ is a
  special metaplectic
  operator. If we introduce the symplectic matrix
 \begin{equation}\label{mJ}
J=\bmatrix 0&I_d\\-I_d&0\endbmatrix,
\end{equation}
then $\cF$ is (up to a phase factor)
the unitary metaplectic operator
corresponding to $J$,
$$
\mu\left(J\right)=(-i)^{d/2}\mc{F}. $$

  \par
  A fundamental tool to
  achieve these estimates is
  represented by the analysis
  of the dilation operator
  $f(x)\mapsto f(Ax)$, for a real invertible $d\times d$ matrix $A\in
  GL(d,\R)$,
  with bounds on its norm in terms of
  spectral invariants of  $A$. In the framework of modulation spaces
  such an investigation  was
  recently developed  in
  the scalar case $A=\lambda
  I$ by
  Sugimoto and Tomita
  \cite{sugimototomita, sugitomita2}. In Section 3 we  study this
  problem for a general matrix
  $A \in GL(d,\R)$  for both  modulation and Wiener amalgam spaces.
  In particular,   we extend
  the results in \cite{sugimototomita} to the case of a
  symmetric matrix $A$.\par
\vskip0.2truecm
  In the second part of the
  paper we present some
  natural
  applications to partial
  differential equations {\it
  with variable
  coefficients}. Precisely,
  we study the Cauchy problem
  for the Schr\"odinger
  equation with a quadratic
  Hamiltonian, namely
  \begin{equation}\label{C1}
\begin{cases} i \displaystyle\frac{\partial
u}{\partial t} +H u=0\\
u(0,x)=u_0(x),
\end{cases}
\end{equation}
where $H$ is the Weyl quantization of a quadratic form on
$\R^d\times \R^d$. The most interesting case is certainly the
Schr\"odinger equation with a quadratic potential. Indeed, the
solution $u(t,x)$ to  \eqref{C1} is given by $$u(t,x)=e^{itH}u_0, $$
where the operator $e^{itH}$ is  a metaplectic operator, so that
 the estimates resulting from the previous sections provide at once
fixed-time estimates for the solution
$u(t,x)$, in terms of the initial datum
$u_0$. An example is provided by the
Harmonic Oscillator
$H=-\frac{1}{4\pi}\Delta+\pi|x|^2$
(see, e.g.,
\cite{folland,helffer,koch}), for which
we deduce the dispersive estimate
\begin{equation}\label{ddd2}
\|e^{itH} u_0\|_{\fui}\lesssim|\sin t|^{-d}\|u_0\|_{\fiu}.
\end{equation}
Another  Hamiltonian we take into account is
$H=-\frac{1}{4\pi}\Delta-\pi|x|^2$ (see \cite{berezinshubin}). In
this case, we show
\begin{equation}\label{d10}
\|e^{itH_\A}
u_0\|_{\fui}\lesssim\left(\frac{1+|\sinh
t|}{\sinh^2 t}\right)^{\frac{d}{2}}
\|u_0\|_{\fiu}.
\end{equation}
In Section $5$ we shall combine these
estimates with orthogonality arguments
as in \cite{cordero,keel} to obtain
space-time estimates: the so-called
Strichartz estimates (for the classical
theory in Lebesgue spaces, see
\cite{GinibreVelo92,Kato70,keel,Yajima87}).
For instance, the homogeneous
Strichartz estimates  achieved for the
Harmonic Oscillator
$H=-\frac{1}{4\pi}\Delta+\pi|x|^2$ read
\[
\|e^{itH} u_0\|_{L^{{q}/{2}}([0,T])
W(\Fur L^{r^\prime},L^r)_x}\lesssim
\|u_0\|_{L^2_x},
\]
for every $T>0$, $4<q,\tilde{q}\leq\infty$, $2\leq
r,\tilde{r}\leq\infty$, such that $2/q+d/r=d/2$, and, similarly, for
$\tilde{q},\tilde{r}$. In the endpoint case $(q,r)=(4, 2d/(d-1))$,
$d>1$, we prove the same estimate with $\Fur L^{r^\prime}$ replaced
by the slightly larger $\Fur L^{r^\prime,2}$, where $L^{r^\prime,2}$
is a Lorentz space (Theorem \ref{prima}).

 The case of the Hamiltonian
$H=-\frac{1}{4\pi}\Delta-\pi|x|^2$ will be detailed in Subsection
$5.2$. Finally, we shall compare all these estimates with the
classical ones in the Lebesgue spaces (Subsection $5.3$).\par

Our analysis  combines techniques from time-frequency analysis
(e.g., convolution relations, embeddings and duality properties of
Wiener amalgam and modulation spaces) with methods from classical
harmonic analysis and PDE's theory (interpolation results,
H\"older-type inequalities, fractional integration theory).

This study carries on the one in \cite{cordero}, developed for the
usual Schr\"odinger equation ($H=-\Delta$).

We record that hybrid spaces like the Wiener amalgam ones had
appeared before as a technical tool in PDE's (see, e.g., Tao
\cite{tao3}). Notice that fixed-time estimates between modulation
spaces in the case $H=-\Delta$ were also considered in
\cite{baoxiang,benyi,benyi2}.\par Finally we observe that, by
combining the Strichartz estimates in the present paper with
arguments from functional analysis as in \cite{Dancona05},
wellposedness for Schr\"odinger equations as above can also be
proved, with an additional potential term in suitable Wiener amalgam
spaces (see also \cite[Section 6]{cordero}). However, here we do not
give details on this subject, that will be studied in a subsequent
paper.

\vskip1truecm

\textbf{Notation.} We define
$|x|^2=x\cdot x$, for
$x\in\Ren$, where $x\cdot
y=xy$ is the scalar product
on $\Ren$. The space of
smooth functions with compact
support is denoted by
$\cC_0^\infty(\rd)$, the
Schwartz class is
$\sch(\Ren)$, the space of
tempered distributions
$\sch'(\Ren)$.    The Fourier
transform is normalized to be
${\hat
  {f}}(\o)=\Fur f(\o)=\int
f(t)e^{-2\pi i t\o}dt$.
 Translation and modulation operators ({\it time and frequency shifts}) are defined, respectively, by
$$ T_xf(t)=f(t-x)\quad{\rm and}\quad M_{\o}f(t)= e^{2\pi i \o
 t}f(t).$$
We have the formulas
$(T_xf)\hat{} = M_{-x}{\hat
{f}}$, $(M_{\o}f)\hat{}
=T_{\o}{\hat {f}}$, and
$M_{\o}T_x=e^{2\pi i
x\o}T_xM_{\o}$. The notation
$A\lesssim B$ means $A\leq c
B$ for a suitable constant
$c>0$, whereas $A \asymp B$
means $c\inv A \leq B \leq c
A$, for some $c\geq 1$. The
symbol $B_1 \hookrightarrow
B_2$ denotes the continuous
embedding of the linear space
$B_1$ into $B_2$.
\section {Function spaces and preliminaries}
\subsection{Lorentz spaces}
 (\cite{stein93,steinweiss}). We
recall that the Lorentz space
$L^{p,q}$ on $\mathbb{R}^d$
is defined as the space of
temperate distributions $f$
such that
\[
\|f\|^\ast_{pq}=\left(\frac{q}{p}\int_0^\infty
[t^{1/p} f^\ast
(t)]^q\frac{dt}{t}\right)^{1/q}<\infty,
\]
when $1\leq p<\infty$, $1\leq
q<\infty$, and
\[
\|f\|^\ast_{pq}=\sup_{t>0}
t^{1/p}f^\ast(t)<\infty
\]
when $1\leq p\leq\infty$, $q=\infty$. Here, as usual,
$\lambda(s)=|\{|f|>s\}|$ denotes the distribution function of $f$
and $f^\ast(t)=\inf\{s:\lambda(s)\leq t\}$.\par One has
$L^{p,q_1}\hookrightarrow L^{p,q_2}$ if $q_1\leq q_2$, and
$L^{p,p}=L^p$. Moreover, for $1<p<\infty$ and $1\leq q\leq\infty$,
$L^{p,q}$ is a normed space and its norm $\|\cdot\|_{L^{p,q}}$ is
equivalent to the above quasi-norm $\|\cdot\|^\ast_{pq}$. \par We
now recall the following generalized Hardy-Littlewood-Sobolev
fractional integration theorem (see e.g. \cite[page 119]{stein} and
\cite[Theorem 2, page 139]{triebel}), which will be used in the
sequel (the original fractional integration theorem corresponds to
the model case of convolution by $K(x)=|x|^{-\alpha}\in
L^{d/\alpha,\infty}$, $0<\a<d$).
\begin{proposition}\label{convlor}
Let $1\leq p<q<\infty$,
$0<\alpha<d$, with
$1/p=1/q+1-\alpha/d$. Then,
\begin{equation}\label{conv1}
L^p(\rd)\ast
L^{d/\alpha,\infty}(\rd)\hookrightarrow
L^q(\rd).
\end{equation}
\end{proposition}
\subsection{Wiener amalgam
spaces}
(\cite{feichtinger80,feichtinger83,
feichtinger90,fournier-stewart85,Fei98}).
Let $g \in \cC_0^\infty$ be a
test function that satisfies
$\|g\|_{L^2}=1$. We will
refer to $g$ as a window
function. Let $B$ one of the
following Banach spaces:
$L^p, \cF L^p$, $1\leq p\leq
\infty$,
 $L^{p,q}$, $1<p<\infty$, $1\leq q\leq \infty$,
  possibly valued in a Banach space, or also
  spaces obtained from these by real or
  complex interpolation.
Let $C$ be one of the following Banach spaces: $L^p$, $1\leq
p\leq\infty$, or $L^{p,q}$,
 $1<p<\infty$, $1\leq q\leq \infty$, scalar-valued.
For any given function $f$ which is locally in $B$ (i.e. $g f\in B$,
$\forall g\in\cC_0^\infty$), we set $f_B(x)=\| fT_x g\|_B$.

The {\it Wiener amalgam space} $W(B,C)$ with local component $B$ and
global component  $C$ is defined as the space of all functions $f$
locally in $B$ such that $f_B\in C$. Endowed with the norm
$\|f\|_{W(B,C)}=\|f_B\|_C$, $W(B,C)$ is a Banach space. Moreover,
different choices of $g\in \cC_0^\infty$  generate the same space
and yield equivalent norms.

If  $B=\Fur L^1$ (the Fourier algebra),  the space of admissible
windows for the Wiener amalgam spaces $W(\Fur L^1,C)$ can be
enlarged to the so-called Feichtinger algebra $W(\Fur L^1,L^1)$.
Recall  that the Schwartz class $\sch$
  is dense in $W(\Fur L^1,L^1)$.\par
We use the following definition of mixed Wiener amalgam norms. Given
a measurable function $F$ of the two variables $(t,x)$ we set
\[
\|F\|_{W(L^{q_1},L^{q_2})_tW(\Fur L^{r_1},L^{r_2})_x}= \|
\|F(t,\cdot)\|_{W(\Fur L^{r_1},L^{r_2})_x}\|_{W(L^{q_1},L^{q_2})_t}.
\]
Observe  that \cite{cordero}
\[
\|F\|_{W(L^{q_1},L^{q_2})_tW(\Fur
L^{r_1},L^{r_2})_x}=
\|F\|_{W\left(L^{q_1}_t(W(\Fur
L^{r_1}_x,L^{r_2}_x)),L^{q_2}_t\right)}.
\]
The following properties of Wiener amalgam spaces  will be
frequently used in the sequel.
\begin{lemma}\label{WA}
  Let $B_i$, $C_i$, $i=1,2,3$, be Banach spaces  such that $W(B_i,C_i)$ are well
  defined. Then,
  \begin{itemize}
  \item[(i)] \emph{Convolution.}
  If $B_1\ast B_2\hookrightarrow B_3$ and $C_1\ast
  C_2\hookrightarrow C_3$, we have
  \begin{equation}\label{conv0}
  W(B_1,C_1)\ast W(B_2,C_2)\hookrightarrow W(B_3,C_3).
  \end{equation}
   In particular, for every
$1\leq p, q\leq\infty$, we have
\begin{equation}\label{p2}
\|f\ast u\|_{\fpq}\leq\|f\|_{W(\Fur L^\infty,L^1)}\|u\|_{\fpq}.
\end{equation}
  \item[(ii)]\emph{Inclusions.} If $B_1 \hookrightarrow B_2$ and $C_1 \hookrightarrow C_2$,
   \begin{equation*}
   W(B_1,C_1)\hookrightarrow W(B_2,C_2).
  \end{equation*}
  \noindent Moreover, the inclusion of $B_1$ into $B_2$ need only hold ``locally'' and the inclusion of $C_1 $ into $C_2$  ``globally''.
   In particular, for $1\leq p_i,q_i\leq\infty$, $i=1,2$, we have
  \begin{equation}\label{lp}
  p_1\geq p_2\,\mbox{and}\,\, q_1\leq q_2\,\Longrightarrow W(L^{p_1},L^{q_1})\hookrightarrow
  W(L^{p_2},L^{q_2}).
  \end{equation}
  \item[(iii)]\emph{Complex interpolation.} For $0<\theta<1$, we
  have
\[
  [W(B_1,C_1),W(B_2,C_2)]_{[\theta]}=W\left([B_1,B_2]_{[\theta]},[C_1,C_2]_{[\theta]}\right),
  \]
if $C_1$ or $C_2$ has
absolutely continuous norm.
    \item[(iv)] \emph{Duality.}
    If $B',C'$ are the topological dual spaces of the Banach spaces $B,C$ respectively, and
    the space of test functions $\cC_0^\infty$ is dense in both $B$ and $C$, then
\begin{equation}\label{duality}
W(B,C)'=W(B',C').
\end{equation}
  \end{itemize}
  \end{lemma}
\begin{proposition}\label{p1}
For every $1\leq p\leq q\leq \infty$, the Fourier transform $\Fur $
maps $\fqp$ in $\fpq$ continuously.
\end{proposition}
\noindent The proof of all these results can be found in
  (\cite{feichtinger80,feichtinger83,feichtinger90,Heil03}).\par
The subsequent result of real interpolation is proved in
\cite{cordero}.
\begin{proposition}\label{inter9}
Given two local components
$B_0,B_1$ as above, for every
$1\leq p_0,p_1<\infty$,
$0<\theta<1$,
$1/p=(1-\theta)/p_0+\theta/p_1$,
and $p\leq q$ we have
\[
W\left((B_0,B_1)_{\theta,q},L^p\right)\hookrightarrow\left(W(B_0,L^{p_0}),W(B_1,L^{p_1})\right)_{\theta,q}.
\]
\end{proposition}

\subsection{Modulation
spaces}(\cite{grochenig}).
Let $g\in\cS$ be a non-zero
window function. The
short-time Fourier transform
(STFT) $V_gf$ of a
function/tempered
distribution $f$ with respect
to the the window $g$ is
defined by
\[
V_g f(z,\o)=\int e^{-2\pi i \o y}f(y)g(y-z)\,dy,
\]
i.e.,  the  Fourier transform $\cF$ applied to $fT_zg$.\par For
$1\leq p, q\leq\infty$, the modulation space $M^{p,q}(\R^n)$ is
defined as the space of measurable functions $f$ on $\R^n$ such that
the norm
\[
\|f\|_{M^{p,q}}=\|\|V_gf(\cdot,\o) \|_{L^p}\|_{L^q_\omega}
\] is
finite. Among the properties of \modsp s, we record that
$M^{2,2}=L^2$, $M^{p_1,q_1}\hookrightarrow M^{p_2,q_2}$, if $p_1\leq
p_2$ and $q_1\leq q_2$. If $p,q<\infty$, then
$(M^{p,q})'=M^{p',q'}$.\par For comparison, notice that the norm in
the Wiener amalgam spaces $W(\Fur L^p,L^q)$ reads
\[
\|f\|_{W(\Fur L^p,L^q)}=\|\|V_gf(z,\cdot) \|_{L^p}\|_{L^q_z}.
\]
The relationship between modulation and Wiener amalgam spaces is
expressed by  the following result.
\begin{proposition}\label{vm} The
Fourier transform establishes
an isomorphism  $\Fur:
M^{p,q}\to W(\Fur
L^p,L^q)$.\par
\end{proposition}
Consequently, convolution properties of modulation spaces can be
translated into point-wise multiplication properties of Wiener
amalgam spaces, as shown below.
\begin{proposition}\label{p3}
For every $1\leq
p,q\leq\infty$ we have
\[
\|fu\|_{\fpq}\leq
\|f\|_{W(\Fur
L^1,L^\infty)}\|u\|_{\fpq}.
\]
\end{proposition}
\begin{proof}
From Proposition \ref{vm}, the estimate to prove is equivalent to
\[
\|\hat{f}\ast
\hat{u}\|_{M^{p,q}}\leq
\|\hat{f}\|_{M^{1,\infty}}\|\hat{u}
\|_{M^{p,q}},
\]
but this a special case  of \cite[Proposition 2.4]{CG02}.
\end{proof}\par\noindent
The  characterization of the $M^{2,\infty}$-norm in \cite[Lemma
3.3]{sugimototomita} can be rephrased  in our context as follows.
\begin{lemma}Suppose that
$\f\in\mathcal{S}(\R^d)$ is a real-valued function satisfying
$\f=1$ on $[-1/2,1/2]^d$, ${\rm supp}\,\f\subset[-1,1]^d$,
$\f(t)=\phi(-t)$ and $\sum_{k\in\mathbb{Z}^d} \f(t-k)=1$ for all
$t\in\R^d$.  Then
\begin{equation}\label{m1}
\|f\|_{M^{2,\infty}}\asymp\sup_{k\in \mathbb{Z}^d}\|(M_k\Phi)\ast
f\|_{L^2},
\end{equation}
for all $f\in M^{2,\infty}$, where $\Phi=\Fur^{-1}\f$.
\end{lemma}

To compute the $\mpq$-norm we shall often use the \emph{duality}
technique, justified by the result below (see \cite[Proposition
11.3.4 and Theorem 11.3.6]{grochenig} and  \cite[Relation
$(2.1)$]{sugimototomita}).
\begin{lemma} Let $\f\in\cS(\rd)$, with $\|\f\|_2=1$, $1\leq p,q<\infty.$ Then
$(M^{p,q})^*=M^{p',q'}$, under the duality
\begin{equation}\label{duality}
\langle f,g\rangle=\langle V_\f f, V_\f g\rangle=\int_{\R^{2d}} V_\f
f(x,\omega) \overline{V_\f g(x,\omega)}\,dx\,d\o,
\end{equation}
for $f\in M^{p,q}$, $g\in
M^{p',q'}$.
\end{lemma}
\begin{lemma}
Assume $1<p,q\leq\infty$ and
$f\in M^{p,q}$. Then
\begin{equation}\label{du}
\|f\|_{M^{p,q}}=\sup_{\|g\|_{M^{p',q'}}\leq 1}|\langle f,g\rangle|.
\end{equation}
\end{lemma}
Notice that  \eqref{du} still holds true whenever $p=1$ or $q=1$ and
$f\in\mathcal{S}(\R^d)$, simply by extending  \cite[Theorem
3.2.1]{grochenig} to the duality ${}_{\cS'}\la\cdot,\cdot\ra_{\cS}$.
\par
Finally we recall the behaviour of
modulation spaces with respect to
complex interpolation (see
\cite[Corollary 2.3]{feichtinger83}.
\begin{proposition}\label{cintm}
Let $1\leq p_1,p_2,q_1,q_2\leq\infty$,
with $q_2<\infty$. If $T$ is a linear
operator such that, for $i=1,2$,
\[
\|Tf\|_{M^{p_i,q_i}}\leq
A_i\|f\|_{M^{p_i,q_i}}\quad \forall
f\in M^{p_i,q_i},
\]
then
\[
\|Tf\|_{M^{p,q}}\leq
CA_1^{1-\theta}A_2^\theta\|f\|_{M^{p,q}}\quad
\forall f\in M^{p,q},
\]
where $1/p=(1-\theta)/p_1+\theta/p_2$,
$1/q=(1-\theta)/q_1+\theta/q_2$, $0<
\theta<1$ and $C$ is independent of
$T$.

\end{proposition}
\subsection{The metaplectic
representation}
(\cite{folland}). The
symplectic group is defined
by
$$
\Spnr=\left\{g\in\Gltwonr:\;^t\!gJg=J\right\},
$$
where the symplectic matrix
$J$ is defined in \eqref{mJ}.
 The metaplectic or Shale-Weil representation
$\mu$ is a unitary
representation of the (double
cover of the) symplectic
group $\Spnr$ on $\lrd$. For
elements of $\Spnr$ in
special form, the metaplectic
representation can be
computed explicitly. For
$f\in L^2(\R^d)$, we have
\begin{align}
\mu\left(\bmatrix A&0\\
0&\;^t\!A^{-1}\endbmatrix\right)f(x)
&=(\det A)^{-1/2}f(A^{-1}x)\nonumber\\
\mu\left(\bmatrix I&0\\
C&I\endbmatrix\right)f(x)
&=\pm e^{i\pi\langle
Cx,x\rangle}f(x)\label{lower}.
\end{align}
The symplectic  algebra
$\spdr$ is the set of all
$2d\times 2d$ real matrices
$\A$   such that $e^{t \A}
\in \Spnr$ for all
$t\in\R$.\par The following
formulae for the metaplectic
representation can be found
in  \cite[Theorems 4.51 and
4.53]{folland}.
\begin{proposition} Let
$\mathcal{A}=\begin{pmatrix}
A&B\\C&D\end{pmatrix}\in
Sp(d,\R)$.

\noindent  (i) If  $\det
B\not=0$ then
\begin{equation}\label{f3}
\mu(\A)f(x)=i^{d/2}(\det
B)^{-1/2}\int e^{-\pi i
x\cdot DB^{-1} x+2\pi i
y\cdot B^{-1}x-\pi i y\cdot
B^{-1}A y} f(y)\,dy.
\end{equation}
\noindent (ii) If  $\det A\not=0$,
\begin{equation}\label{f4}
\mu(\A)f(x)=(\det
A)^{-1/2}\int e^{-\pi i
x\cdot CA^{-1}x+2\pi i
\xi\cdot A^{-1} x+\pi
i\xi\cdot
A^{-1}B\xi}\hat{f}(\xi)\,d\xi.
\end{equation}
\end{proposition}
The following hybrid formula will be also used in the sequel.
\begin{proposition}
If $\mathcal{A}=\begin{pmatrix} A&B\\C&D\end{pmatrix}\in Sp(d,\R)$,
$\det B\not=0$ and $\det A\not=0$, then
\begin{equation}\label{f5}
\mu(\A) f(x)=(-i\det
B)^{-1/2} e^{-\pi ix\cdot
CA^{-1}x}\left(e^{-\pi i
y\cdot B^{-1}A y}\ast
f\right)(A^{-1}x).
\end{equation}
\end{proposition}
\begin{proof}
By \eqref{f4} we can write
\begin{align}
\mu(\A) f(x)&=(\det A)^{-1/2}
e^{-\pi ix\cdot CA^{-1}x}\int
e^{2\pi i\xi\cdot A^{-1}x}
\Fur\left(\Fur^{-1}e^{\pi
i\xi\cdot
A^{-1}B\xi}\right)\hat{f}(\xi)\,d\xi\nonumber\\
&=(-i\det B)^{-1/2}e^{-\pi
ix\cdot CA^{-1}x}\int e^{2\pi
i\xi\cdot A^{-1}x} \Fur\left(
e^{-\pi iy\cdot B^{-1}A
y}\ast
f\right)(\xi)\,d\xi,\nonumber
\end{align}
where we used the formula
(see \cite[Theorem 2, page
257] {folland})
\[
\Fur^{-1}
\left(e^{i\pi\xi\cdot
A^{-1}B\xi}\right)(y)= (-i\det
A^{-1}B)^{-1/2} e^{-\pi
iy\cdot B^{-1}A y}.
\]
Hence, from the Fourier
inversion formula we obtain
\eqref{f5}.
\end{proof}

\section{ Dilation
 of  Modulation and Wiener Amalgam Spaces} Given a function
$f$ on $\rd$ and $A\in
GL(d,\R)$, we set
$f_A(t)=f(At)$. We also
consider the unitary operator
$\cU_A$ on $\lrd$ defined by
\begin{equation}\label{dilop}
    \cU_A f(t)=|\det A|^{1/2} f(A t)=|\det A|^{1/2} f_A( t).
\end{equation}
In this section we study the boundedness of this operator on
modulation and Wiener amalgam spaces. We need the following three
lemmata.
\begin{lemma}\label{l1} Let $A\in GL(d,\R)$, $\f(t)=e^{-\pi |t|^2}$, then
\begin{equation*}
V_{\f}\f_A \phas= (\det(A^* A+I))^{-1/2} e^{-\pi (I-(A^*
A+I)^{-1})x\cdot x} M_{-((A^* A+I)^{-1})x}e^{-\pi (A^*
A+I)^{-1}\o\cdot\o}.
\end{equation*}
\end{lemma}
\begin{proof} By definition of the STFT,
\begin{eqnarray*}
  V_{\f}\f_A \phas &=& \intrd e^{-\pi Ay\cdot y} e^{-2\pi i \o\cdot y}e^{-\pi (y-x)^2}\, dy \\
   &=& e^{-\pi |x|^2}\intrd e^{-\pi (A^* A+I)y\cdot y+2 \pi x\cdot y} e^{-2\pi i \o\cdot y}\, dy.  \\
\end{eqnarray*}
Now, we rewrite the generalized Gaussian above using the
translation and dilation operators, that is
$$e^{-\pi (A^* A+I)y\cdot y+2 \pi x\cdot y}=
(\det(A^*A+I))^{-1/4}(T_{(A^*A+I)^{-1}x}\cU_{(A^*A+I)^{1/2}})\f(y)
$$
and use the properties $\cF \cU_{B}=\cU_{(B^*)^{-1}}\cF$, for
every $B\in GL(d,\R)$ and $\cF T_x= M_{-x}\cF$. Thereby,
\begin{eqnarray*}
  V_{\f}\f_A \phas &=& e^{-\pi( I-(A^*A+I)^{-1})x\cdot x}(\det(A^*A+I))^{-1/4}\cF(T_{(A^*A+I)^{-1}x}\cU_{(A^*A+I)^{1/2}}\f(\o)  \\
   &=&  e^{-\pi( I-(A^*A+I)^{-1})x\cdot x}(\det(A^*A+I))^{-1/2}M_{-(A^*A+I)^{-1}x}e^{-\pi(A^*A+I)^{-1}\o\cdot \o},  \\
\end{eqnarray*}
as desired.
\end{proof}

The result below generalizes \cite[Lemma 1.8]{Toft04}, recaptured
in the special case $A=\lambda I$, $\lambda>0$.
\begin{lemma}\label{l2} Let $1\leq p,q\leq \infty$, $A\in GL(d,\R)$ and  $\f(t)=e^{-\pi |t|^2}$. Then,
\begin{equation}\label{e2}
\|\f_A\|_{\mpq}=p^{-d/(2p)}q^{-d/(2q)}|\det
A|^{-1/p}(\det(A^*A+I))^{-(1-1/q-1/p)/2}.
\end{equation}
\end{lemma}
\begin{proof} Since the modulation space norm is independent  of the choice of the window function, we choose the Gaussian $\f$, so that $\|\f_A\|_{\mpq}\asymp\|V_\f\f_A\|_{L^{p,q}}$. Since
\begin{eqnarray*}
  \intrd e^{-\pi p(I-(A^*A+I)^{-1})x\cdot
x}\,dx &=& (\det(I-(A^*A+I)^{-1})^{-1/2} p^{-d/2} \\
  &=&p^{-d/2}|\det A|^{-1}(\det(A^*A+I))^{1/2}
\end{eqnarray*}
 and, analogously, $\intrd e^{-\pi q(A^*A+I)^{-1}\o\cdot
\o}\,d\o=(\det(A^*A+I))^{1/2} q^{-d/2}$, the result immediately
follows from Lemma \ref{l1}.
\end{proof}\par\noindent
We record \cite[Lemma 11.3.3]{grochenig}:
\begin{lemma}\label{l3} Let $f\in\cS'(\rd)$ and $\f,\psi,\gamma\in\cS(\rd)$. Then,
$$|V_\f f\phas|\leq\frac1{\la \gamma,\psi\ra}(|V_{\psi}f|\ast|V_{\f}\gamma|)\phas\quad \forall \phas\in\rdd.
$$
\end{lemma}
The results above are the
ingredients for the first
dilation property of \modsp s
we are going to present.
\begin{proposition}\label{dilmod} Let $1\leq p,q\leq\infty$ and  $A\in GL(d,\R)$. Then, for every $f\in \mpq(\rd)$,
\begin{equation}\label{pp1}\|f_A\|_{\mpq}\lesssim
|\det
A|^{-(1/p-1/q+1)}(\det(I+A^\ast
A))^{1/2}\|f\|_{\mpq}.
\end{equation}
\end{proposition}
\begin{proof} The proof follows the guidelines  of \cite[Lemma 3.1]{sugimototomita}. First, by
a change of variable,  the dilation is transferred from the function $f$ to the window $\f$:
$$V_\f f_A\phas=|\det A|^{-1} V_{\f_{A^{-1}}} f(A x,(A^*)^{-1}\o).$$
Whence,  performing the change of variables $Ax=u$, $(A^*)^{-1}\o=v$,
\begin{align*}
\|f_A\|_{\mpq}&=|\det A|^{-1}\left(\intrd\left(\intrd | V_{\f_{A^{-1}}}
f(A x,(A^*)^{-1}\o)|^p\,dx\right)^{q/p}d\o\right)^{1/q} \\
&=|\det A|^{-(1/p-1/q+1)}\|V_{\f_{A^{-1}}}\|_{L^{p,q}}.
\end{align*}
Now, Lemma \ref{l3}, written for $\psi(t)=\gamma(t)= \f(t)=e^{-\pi t^2}$, yields  the following majorization
$$|V_{\f_{A^{-1}}}f\phas|\leq \|\f\|_{L^2}^{-2}  (|V_{\f} f|\ast|V_{\f_{A^{-1}}} \f|)\phas.
$$
Finally, Young's Inequality and Lemma \ref{l2} provide the desired result:
\begin{align*}
\|f_A\|_{\mpq}&\lesssim |\det A|^{-(1/p-1/q+1)} \||V_{\f} f|\ast|V_{\f_{A^{-1}}} \f|\|_{L^{p,q}}\\
&\lesssim  |\det A|^{-(1/p-1/q+1)}\|V_{\f} f\|_{L^{p,q}}\|V_{\f_{A^{-1}}} \f\|_{L^{1}}\\
&\asymp  |\det
A|^{-(1/p-1/q+1)}(\det(I+A^*A))^{1/2}\|f\|_{\mpq}.
\end{align*}
\end{proof}

Proposition \ref{dilmod}
generalizes \cite[Lemma
3.1]{sugimototomita}, that
can be recaptured by choosing
the matrix  $A=\lambda I$,
$\lambda>0$.
\begin{corollary}\label{c1}
Let $1\leq p,q\leq\infty$ and  $A\in
GL(d,\R)$. Then, for every $f\in
\fpq(\rd)$,\begin{equation}\label{dilAW0}
\|f_A\|_{\fpq}\lesssim |\det
A|^{(1/p-1/q-1)}(\det(I+A^*
A))^{1/2}\|f\|_{\fpq}.
\end{equation}
\end{corollary}
\begin{proof} It follows immediately from the relation between Wiener amalgam spaces and modulation spaces given by $\fpq=\cF\mpq$ and by the relation $\widehat{(f_A)}=|\det A|^{-1}(\hat{f})_{(A^*)^{-1}}$.
\end{proof}

In what follows we give a
more precise result about the
behaviour of the operator
norm
$\|D_A\|_{\mpq\rightarrow
\mpq}$ in terms of $A$, when
$A$ is a symmetric matrix,
extending the diagonal case
$A=\lambda I$, $\lambda>0$
treated in
\cite{sugimototomita}.  We
shall use the  set and index
terminology of the paper
above. Namely, for $1\leq
p\leq\infty$, let $p'$ be the
conjugate exponent of $p$
($1/p+1/p'=1$). For
$(1/p,1/q)\in [0,1]\times
[0,1]$, we define the subsets
$$ I_1=\max (1/p,1/p')\leq 1/q,\quad\quad I_1^*=\min (1/p,1/p')\geq 1/q,
$$
$$ I_2=\max (1/q,1/2)\leq 1/p',\quad\quad I_2^*=\min (1/q,1/2)\geq  1/p',
$$
$$ I_3=\max (1/q,1/2)\leq 1/p,\quad\quad I_3^*=\min (1/q,1/2)\geq
1/p,
$$
as shown in Figure 1:
\vspace{1.2cm}
 \begin{center}
           \includegraphics{figSchr1.1}
            \\
           $ $
\end{center}
 \begin{center}{\quad \quad\quad\quad \quad $0<|\lambda|\leq 1$\hfill   $|\lambda|\geq 1$}\quad\quad\quad\quad\quad\quad\quad\quad\quad
           \end{center}
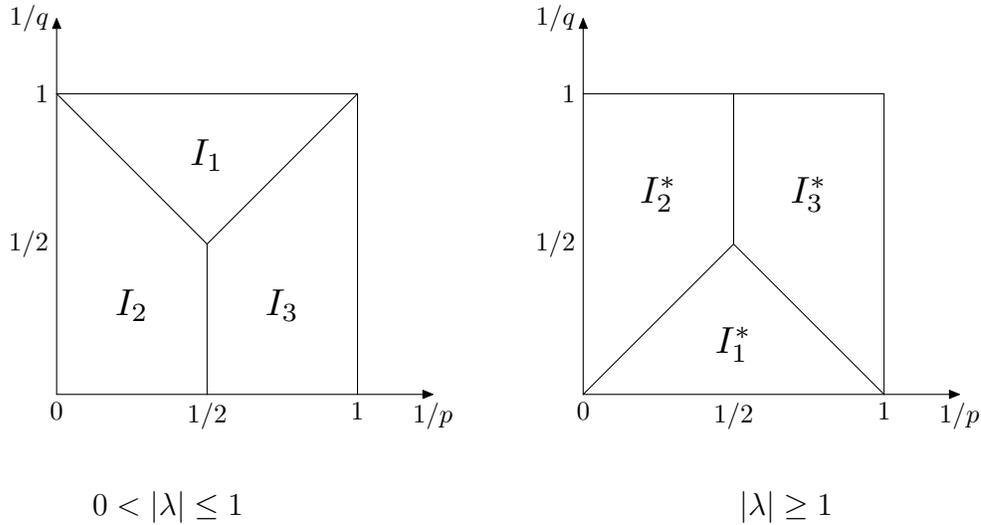
           \begin{center}{ Figure 1. The index sets. }
           \end{center}
  \vspace{1.2cm}

We introduce the indices:
$$ \mu_1(p,q)=\begin{cases}-1/p &  \quad {\mbox{if}} \quad (1/p,1/q)\in  I_1^*,\\
 1/q-1 &   \quad {\mbox{if}}  \quad (1/p,1/q)\in  I_2^*,\\
 -2/p +1/q&  \quad  {\mbox{if}}  \quad (1/p,1/q)\in  I_3^*,\\
 \end{cases}
 $$
and
$$ \mu_2(p,q)=\begin{cases}-1/p &  \quad {\mbox{if}} \quad (1/p,1/q)\in  I_1,\\
 1/q-1 &   \quad {\mbox{if}}  \quad (1/p,1/q)\in  I_2,\\
 -2/p +1/q&  \quad  {\mbox{if}}  \quad (1/p,1/q)\in  I_3.\\
 \end{cases}
 $$
The above mentioned result by
\cite{sugimototomita} reads as follows:

\begin{theorem}\label{dilprop}
Let $1\leq p,q\leq\infty,$ and
$A=\lambda I$, $\lambda\not=0$.\par
 (i) We have
$$\| f_A\|_{\mpq}\lesssim |\lambda|^{d\mu_1(p,q)}
\|f\|_{\mpq},\quad\quad\forall\
|\lambda|\geq 1,\ \forall
f\in\mpq(\rd).
$$
Conversely, if there exists  $\a>0$ such that
$$\| f_A\|_{\mpq}\lesssim |\lambda|^{\a}\|f\|_{\mpq},
\quad\quad \forall\ |\lambda|\geq 1,\
\forall f\in\mpq(\rd),
$$
then $\a\geq d\mu_1(p,q)$.\par
 (ii) We
have
$$\| f_A\|_{\mpq}\lesssim |\lambda|^{d\mu_2(p,q)}
\|f\|_{\mpq},\quad\quad\forall\
0<|\lambda|\leq 1,\ \forall
f\in\mpq(\rd).
$$
Conversely, if there exists  $\b>0$ such that
$$\| f_A\|_{\mpq}\lesssim |\lambda|^{\b}\|f\|_{\mpq},
\quad\quad \forall\ 0<|\lambda|\leq 1,\
\forall f\in\mpq(\rd),
$$
then $\b\leq d\mu_2(p,q)$.
\end{theorem}
Here is our extension.
\begin{theorem}\label{t5}
Let $1\leq p,q\leq\infty$. There exists a constant $C>0$ such
that, for every symmetric matrix $A\in { GL}(d,\R)$, with
eigenvalues $\lambda_1,\ldots,\lambda_d$, we have
\begin{equation}\label{dilA}
\| f_A\|_{\mpq}\leq C\prod_{j=1}^d
(\mal)^{\mu_{1}(p,q)}(\mil)^{\mu_{2}(p,q)} \|f\|_{\mpq},
\end{equation}
for every $f\in\mpq(\rd)$.
\par Conversely, if there
exist $\a_j>0,\b_j>0$ such
that, for every
$\lambda_j\not=0$,
$$\| f_A\|_{\mpq}\leq C \prod_{j=1}^d
(\mal)^{\alpha_j}(\mil)^{\beta_j}
\|f\|_{\mpq},\quad\quad
\forall f\in\mpq(\rd),
$$
with $A={\rm diag}[\lambda_1,\ldots,\lambda_d]$, then
$\a_j\geq\mu_1(p,q)$ and $\b_j\leq \mu_2(p,q)$.
\end{theorem}
\begin{proof}
The necessary conditions are an immediate consequence of the one
dimensional case,  already contained in Theorem \ref{dilprop}.
Indeed, it can be seen by taking $f$ as tensor product of functions
of one variable and by leaving free to vary just one eigenvalue, the
remaining eigenvalues being all equal to one. \par Let us come to
the first part of the theorem. It suffices to prove it in the
diagonal case   $A=D={\rm diag}[\lambda_1,\ldots,\lambda_d]$.
Indeed, since $A$ is symmetric, there exists an orthogonal matrix
$T$ such that $A=T^{-1} D T$, and $D$ is a diagonal matrix. On the
other hand, by Proposition \ref{dilmod}, we have
$\|f_A\|_{\mpq}\lesssim\|f_{T^{-1}D}\|_{\mpq}=\|(f_{T^{-1}})_{D}\|_{\mpq}$
and $\|f_{T^{-1}}\|_{\mpq}\lesssim \|f\|_{\mpq}$; hence  the general
case in \eqref{dilA}  follows from the diagonal  case  $A=D$, with
$f$  replaced  by $f_{T^{-1}}$.

From now onward, $A=D={\rm diag}[\lambda_1,\ldots,\lambda_d]$.\par

 If
the theorem holds true for a pair $(p,q)$, with $(1/p,1/q)\in
[0,1]\times[0,1]$, then it is also true for their dual pair
$(p',q')$ (with $f\in\mathcal{S}$ if $p'=1$ or $q'=1$, see
\eqref{du}). Indeed,
\begin{align*}
\|f_D\|_{M^{p',q'}}&=\sup_{\|g\|_{M^{p,q}}\leq 1}|\langle
f_D,g\rangle|= |\det D|^{-1}\sup_{\|g\|_{M^{p,q}}\leq 1}|\langle
f,g_{D^{-1}}\rangle|\\
& \leq |\det D|^{-1}\|f\|_{M^{p',q'}} \sup_{\|g\|_{M^{p,q}}\leq
1}\|g_{D^{-1}}\|_{M^{p,q}}\\&\lesssim
\prod_{j=1}^d|\lambda_j|^{-1}\prod_{j=1}^d(\max\{1,|\lambda_j|^{-1}\})^{
\mu_1(p,q)}(\min\{1,|\lambda_j|^{-1}\})^{\mu_2(p,q)}
\|f\|_{M^{p',q'}}\\
&=\prod_{j=1}^d(\max\{1,|\lambda_j|\})^{
\mu_1(p',q')}(\min\{1,|\lambda_j|\})
^{\mu_2(p',q')}\|f\|_{M^{p',q'}},
\end{align*}
for the index functions $\mu_1$ and $\mu_2$ fulfill:
\begin{equation}\label{indexr}
    \mu_1(p',q')=-1-\mu_2(p,q), \quad \mu_2(p',q')=-1-\mu_1(p,q).
\end{equation}
Hence it suffices to prove
 the  estimate \eqref{dilA} for the
case $p\geq q$. Notice that the estimate in $M^{1,q'}$, $q'>1$, are
proved for Schwartz functions only, but they extend to all functions
in $M^{1,q'}$, $q'<\infty$, for  $\mathcal{S}(\R^d)$ is dense in
$M^{1,q'}$.  The uncovered case $(1,\infty)$ will be verified
directly at the end of the proof.

From Figure $1$  it is clear that  the estimate \eqref{dilA}
 for the points in the upper triangles follows  by
complex interpolation (Proposition
\ref{cintm}) from the diagonal case
$p=q$, and the two cases
$(p,q)=(\infty,1)$ and $(p,q)=(2,1)$.
 \vspace{1.2cm}
 \begin{center}
           \includegraphics{figSchr3.1}
            \\
           $ $
\end{center}
 \begin{center}{\quad \quad\quad\quad \quad \qquad$p\geq q$\hfill \quad  $p\leq q$}\quad\quad\quad\quad\quad\quad\quad\quad\quad\quad\quad
           \end{center}
\vspace{0.2cm}
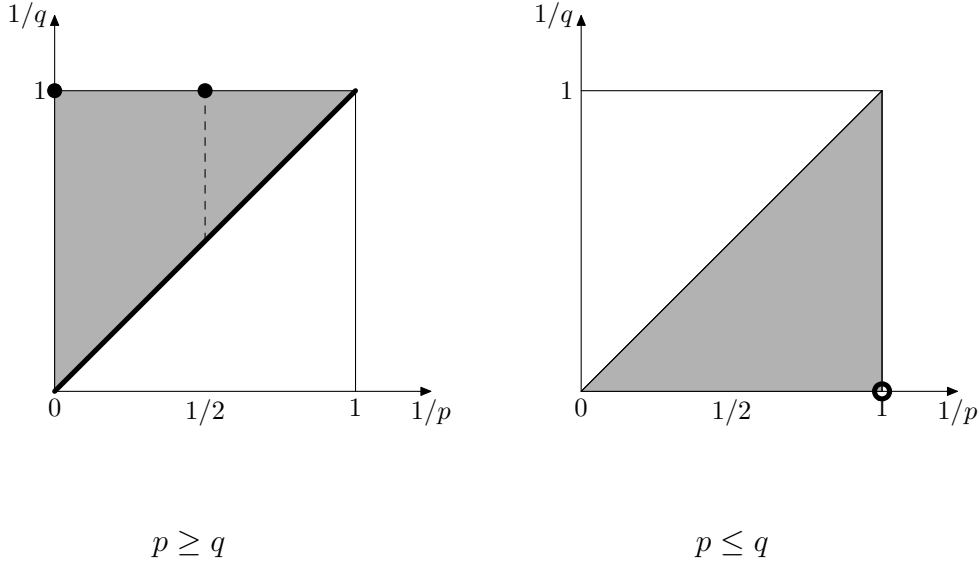
           \begin{center}{ Figure 2. The complex interpolation and the duality method. }
           \end{center}
  \vspace{1.2cm}

\par\medskip\noindent
\emph {Case $p=q$}. If $d=1$ the claim is true by Theorem
\ref{dilprop} in dimension $d= 1$. We then use the induction method.
Namely, we assume that \eqref{dilA} is fulfilled in dimension $d-1$
and prove that still holds in dimension $d$.\\ For $x,\o\in \rd$, we
write $x=(x',x_d)$, $\o=(\o',\o_d)$, with $x',\o'\in \R^{d-1}$,
$x_d,\o_d\in\R$, $D'={\rm diag}[\lambda_1,\ldots,\lambda_{d-1}]$,
and choose the Gaussian
$\f(x)=e^{-\pi|x|^2}=e^{-\pi|x'|^2}e^{-\pi|x_d|^2}=\f'(x')\f_d(x_d)$
as window function. Observe that $V_\f f_D$ admits the two
representations
\begin{align*}
V_\f f_D(x',x_d,\o',\o_d)&=\intrd f(\lambda_1 t_1,\dots,\lambda_d
t_d) \overline{M_{\o'}T_{x'}\f'(t')}
\overline{M_{\o_d}T_{x_d}\f_d(t_d)}\,dt' dt_d\nonumber\\
&= V_{\f'}((F_{x_d,\o_d,\lambda_d})_{D'})\\
&= V_{\f_d}((G_{x',\o',D'})_{\lambda_d})
\end{align*}
where
$$F_{x_d,\o_d,\lambda_d}(t')=V_{\f_d}(f(
t',\lambda_d \cdot)(x_d,\o_d),\quad
G_{x',\o',D'}(t_d)=V_{\f'}(f(D' \cdot,
t_d))(x',\o').$$ By the inductive
hypothesis we have
\begin{align*}
\|f_D\|_{M^{p,p}(\rd)}&=\|V_\f f_D\|_{L^p(\rdd)}\nonumber\\
&=\left(\int_{\R^2}\left(\int_{\R^{2(d-1)}}|V_{\f'}((F_{x_d,\o_d,\lambda_d})_{D'})(x',\o')|^p
dx'd\o'\right)dx_d d\o_d\right)^{1/p}\\ &\lesssim
\prod_{j=1}^{d-1}(\mal)^{\mu_{1}(p,p)} (\mil)^{\mu_{2}(p,p)}\\
&\qquad\qquad\qquad\cdot \left(\intrdd
|V_{\f'}(F_{x_d,\o_d,\lambda_d})
(x',\o')|^p\,dx\,d\o\right)^{1/p}\\
&= \prod_{j=1}^{d-1}(\mal)^{\mu_{1}(p,p)} (\mil)^{\mu_{2}(p,p)}\\
&\qquad\qquad\qquad\cdot \left(\int_{\R^{2(d-1)}} \left(\int_{\R^2}
 |V_{\f_d}((G_{x',\o',I})_{\lambda_d})
(x_d,\o_d)|^p\,dx_d\,d\o_d\right)dx'd\o'\right)^{1/p}\\
&\lesssim \prod_{j=1}^{d}(\mal)^{\mu_{1}(p,p)} (\mil)^{\mu_{2}(p,p)}
\|f\|_{M^{p,p}(\rd)},
\end{align*}
where in the last raw we used Theorem \ref{dilprop} for $d=1$.
\par\medskip\noindent \emph{ Case
$(p,q)=(2,1)$}. First, we prove the
case $(p,q)=(2,\infty)$ and then obtain
the claim by duality as above, since
$\mathcal{S}$ is dense in $M^{2,1}$.
Namely, we want to show that
\[
\|f_D\|_{M^{2,\infty}}\lesssim
\prod_{j=1}^d
\mal)^{-1/2}(\mil)^{-1}\|f\|_{M^{2,\infty}},
\quad \forall f\in
M^{2,\infty}.
\]
The arguments are similar to \cite[Lemma 3.4]{sugimototomita}. We
use the characterization of the $M^{2,\infty}$-norm in \eqref{m1}
\begin{align}
\|f_D\|_{M^{2,\infty}}&\lesssim |\det D|^{-1/2}\sup_{k\in
\mathbb{Z}^d}\|\f(D\cdot-k)\hat{f}\|_{L^2}\nonumber
\\
&=|\det D|^{-1/2}\sup_{k\in
\mathbb{Z}^d}\|\f (D\cdot-k)
\left(\sum_{l\in\mathbb{Z}^d}\f(\cdot-l)
\right)\hat{f}\|_{L^2}\label{m2}.
\end{align}
Observe that
\begin{eqnarray*}
  \left| \f (D t-k) \left(\sum_{l\in\mathbb{Z}^d}\f(t-l)
\right)\hat{f}(t)\right|^2 &\leq& 4^d
\sum_{l\in\mathbb{Z}^d}\left| \f (D t-k) \f(t-l) \hat{f}(t)\right|^2 \\
  &=& 4^d
\sum_{l\in\Lambda_k}\left| \f (D t-k) \f(t-l) \hat{f}(t)\right|^2
\end{eqnarray*}
 where
$$\Lambda_k=\left\{l\in\mathbb{Z}^d:|l_j-\frac{k_j}{\lambda_j}|\leq
1+\frac{1}{|\lambda_j|}\right\} $$ and
$$
\#\Lambda_k\leq C\prod_{j=1}^d \mil^{-1},\quad \forall
k\in\mathbb{Z}^d$$ ($C$ being a constant depending on $d$ only).
Since $|\lambda_j|=\mal\mil$, the expression on the right-hand side
of \eqref{m2} is dominated by
\[
C'\prod_{j=1}^d(\mal)^{-1/2}(\mil)^{-1}
\sup_{m\in\mathbb{Z}^d}\|(M_m\Phi)\ast f\|_{L^2}.
\]
Thereby the norm equivalence \eqref{m1}
gives the desired estimate.
\par\medskip\noindent
\emph {Case  $(p,q)=(\infty,1)$}. We have to prove that
\[
\|f_D\|_{M^{\infty,1}}\lesssim \prod_{j=1}^d
\mal\|f\|_{M^{\infty,1}},\quad \forall f\in M^{\infty,1}.
\]
This estimate immediately
follows from \eqref{pp1},
written for $A=D={\rm
diag}[\lambda_1,\ldots,\lambda_d]$:
\[
\|f_D\|_{M^{\infty,1}}\lesssim \prod_{j=1}^d
(1+\lambda_j^2)^{1/2}\lesssim \prod_{j=1}^d \mal
\|f\|_{M^{\infty,1}}.
\]

\par\medskip\noindent
\emph {Case  $(p,q)=(1,\infty)$}. We are left to prove that
\[
\|f_D\|_{M^{1,\infty}}\lesssim \prod_{j=1}^d
(\mal)^{-1}(\mil)^{-2}\|f\|_{M^{1,\infty}},\quad \forall f\in
M^{1,\infty}.
\]
This is again the estimate  \eqref{pp1}, written for $A=D={\rm
diag}[\lambda_1,\ldots,\lambda_d]$:
\[
\|f_D\|_{M^{1,\infty}}\lesssim
\prod_{j=1}^d|\lambda_j|^{-2}\prod_{j=1}^d \mal
\|f\|_{M^{1,\infty}}.
\]
\end{proof}
\begin{corollary}
\label{c5} Let $1\leq p,q\leq\infty$. There exists a constant $C>0$
such that, for every symmetric matrix $A\in { GL}(d,\R)$, with
eigenvalues $\lambda_1,\ldots,\lambda_d$, we have
\begin{equation}\label{dilAW}
\| f_A\|_{W(\Fur L^p,L^q)}\leq C\prod_{j=1}^d
(\mal)^{\mu_{1}(p',q')}(\mil)^{\mu_{2}(p',q')} \|f\|_{W(\Fur
L^p,L^q)},
\end{equation}
for every $f\in {W(\Fur L^p,L^q)}(\rd)$.
\par Conversely, if there
exist $\a_j>0,\b_j>0$ such that, for every $\lambda_j\not=0$,
$$\| f_A\|_{W(\Fur
L^p,L^q)}\leq C \prod_{j=1}^d (\mal)^{\alpha_j}(\mil)^{\beta_j}
\|f\|_{W(\Fur L^p,L^q)},
$$
for every $f\in {W(\Fur L^p,L^q)}(\rd)$, with $A={\rm
diag}[\lambda_1,\ldots,\lambda_d]$, then $\a_j\geq\mu_1(p',q')$ and
$\b_j\leq \mu_2(p',q')$.\par
\end{corollary}
\begin{proof} It is a mere consequence of Theorem \ref{t5} and
the index relation \eqref{indexr}. Namely,
\begin{eqnarray*}
 \| f_A\|_{W(\Fur
L^p,L^q)} &=&  \|\widehat{ f_A}\|_{\mpq}=|\det A|^{-1}
 \|\hat{ f}_{A^{-1}}\|_{\mpq}\\
   &\leq& C \prod_{j=1}^d|\lambda_j|^{-1} \prod_{j=1}^d(\max\{1,|\lambda_j|^{-1}\})^{\mu_1(p,q)}
   (\min\{1,|\lambda_j|^{-1}\})^{\mu_2(p,q)}\|\hat{ f}\|_{\mpq}\\
   &=&  C \prod_{j=1}^d
(\mal)^{\mu_{1}(p',q')}(\mil)^{\mu_{2}(p',q')} \|f\|_{W(\Fur
L^p,L^q)}.
\end{eqnarray*}
The necessary conditions use the same argument.
\end{proof}
\section{Action of metaplectic operators
on Wiener amalgam spaces} In this section we study the continuity
property of metaplectic operators on Wiener amalgam spaces, giving
bounds on their norms. Here is our first result.
\begin{theorem}\label{t1} Let
$\mathcal{A}=\begin{pmatrix}
A&B\\C&D\end{pmatrix}\in Sp(d,\R)$, and
$1\leq p\leq q\leq\infty$.\par
\noindent (i) If $\det B\not=0$,
 then
\begin{equation}\label{f1}
\|\mu(\mathcal{A})f\|_{\fpq}\lesssim
\a(\mathcal{A},p,q)\|f\|_{\fqp},
\end{equation}
where
\begin{equation}\label{alfa}
\a(\mathcal{A},p,q)= |\det B|^{1/q-1/p-3/2}|\det(I+B^*B)(B+iA)
(B+iD)|^{1/2}.
\end{equation}
 (ii) If $\det A,\det B\not=0$, then
\begin{equation}\label{f2}
\|\mu(\mathcal{A})f\|_{W(\cF
L^1,L^\infty)}\lesssim
\beta(\mathcal{A})\|f\|_{W(\cF
L^\infty,L^1)},
\end{equation}
with
\begin{equation}\label{beta}
\beta(\mathcal{A})= |\det A|^{-3/2}|\det
B|^{-1}|\det(I+A^*A)(B+iA) (A+iC)|^{1/2}.
\end{equation}
\end{theorem}
If the matrices  $A$  or $B$ are symmetric, Theorem \ref{t1} can be
sharpened as follows.
\begin{theorem}\label{ttt1}
Let $\mathcal{A}=\begin{pmatrix}
A&B\\C&D\end{pmatrix}\in Sp(d,\R)$, and
$1\leq p\leq q\leq\infty$.\par
\noindent (i) If $\det B\not=0$,
$B^\ast=B$, with eigenvalues
$\lambda_1,\ldots,\lambda_d$,
 then
\begin{equation}\label{f1z}
\|\mu(\mathcal{A})f\|_{\fpq}\lesssim
\a'(\mathcal{A},p,q)\|f\|_{\fqp},
\end{equation}
where
\begin{align}
\a'(\mathcal{A},p,q)&= |\det(B+iA) (B+iD)|^{1/2}\nonumber\\
&\qquad\qquad\qquad\qquad\cdot\prod_{j=1}^d(\mal)^{\mu_1(p,q)-1/2}(\mil)^{\mu_2(p,q)-1/2}.\label{f1za}
\end{align}
 (ii) If $\det A,\det B\not=0$, and
 $A^\ast=A$ with eigenvalues
$\nu_1,\ldots,\nu_d$,
 then
\begin{equation}\label{f2z}
\|\mu(\mathcal{A})f\|_{W(\cF
L^1,L^\infty)}\lesssim
\beta'(\mathcal{A})\|f\|_{W(\cF
L^\infty,L^1)},
\end{equation}
with
\begin{align}
\beta'(\mathcal{A})&=|\det B|^{-1}|\det(B+iA) (A+iC)|^{1/2}\nonumber\\
&\qquad\qquad\qquad\qquad\cdot\prod_{j=1}^d
(\max\{1,|\nu_j|\})^{-1/2}(\min\{1,|\nu_j|\})^{-3/2}.\label{f2zb}
\end{align}
\end{theorem}

We now prove Theorems \ref{t1} and
\ref{ttt1}. We need the following
preliminary result.
\begin{lemma}
Let $R$ be a $d\times d$ real symmetric matrix, and $f(y)=e^{-\pi iR
y\cdot y}$. Then,
\begin{equation}\label{chirpnorm}
\|f\|_{W(\Fur
L^1,L^\infty)}=|\det(I+iR)|^{1/2}.
\end{equation}
\end{lemma}
\begin{proof}
We first compute the short-time Fourier transform of $f$, with
respect to the window $g(y)=e^{-\pi|y|^2}$. We have
\begin{align}
V_g f(x,\xi)&=\int e^{-2\pi i
y\xi} e^{-i\pi Ry\cdot y}
e^{-\pi|y|^2}
\,dy\nonumber\\
&=e^{-\pi|x|^2}\int e^{-2\pi
i
y\cdot(\xi+ix)-\pi(I+iR)y\cdot
y}\,dy\nonumber\\
&=e^{-\pi|x|^2}(\det(I+iR))^{-1/2}
e^{-\pi(I+iR)^{-1}(\o+ix)
\cdot(\o+ix)}, \nonumber
\end{align}
where we used  \cite[Theorem 1,  page 256]{folland}. Hence
\[
|V_g
f(x,\xi)|=|\det(I+iR)|^{-1/2}e^{-\pi
(I+R^2)^{-1}(\xi+Rx)\cdot(\xi+Rx)},
\]
and, performing the change of
variables
$(I+R^2)^{-1/2}(\xi+Rx)=y$,
with $d\o=|\det(I+R^2)|^{1/2}
dy$, we obtain
\begin{equation}\label{fe}
\intrd V_g f\phas d\o=|\det(I+iR)|^{-1/2}(\det (I+R^2))^{1/2}\intrd
e^{-\pi |y|^2}\,d y=|\det(I+iR)|^{1/2}.
\end{equation}
The last equality follows from $(I+iR)=(I+R^2)(I-iR)^{-1}$, so that
$\det (I+iR)^{-1}= \det(I+R^2)^{-1} \det(I-iR)$. Now, relation
\eqref{chirpnorm} is proved by taking the supremum with the respect
to $x\in\R^d$ in \eqref{fe}.
\end{proof}
\begin{proof}[Proof
of Theorem \ref{t1}] \emph{(i)} We use the expression of
$\mu(\A)f$ in formula \eqref{f3}. The estimates below are obtained
by using (in order): Proposition \ref{p3}  with Lemma 4.1, the
estimate \eqref{dilAW0}, Proposition \ref{p1}, and, finally,
Proposition \ref{p3} combined with Lemma 4.1 again:
\begin{align*}
    \|\mu(\A)f\|_{\fpq}    &=|\det B|^{-1/2}\left\|e^{-\pi i x\cdot DB^{-1}x} \cF^{-1}\left(
     e^{-\pi i  y\cdot  B^{-1}A y} f\right)(B^{-1} x)\right\|_{\fpq}\\
&\leq |\det B|^{-1/2}\|e^{-\pi i x\cdot DB^{-1}x}\|_{W(\cF
L^1,L^\infty)}\\
&\qquad\qquad\quad\quad\quad\cdot\left\|\left( \cF^{-1}\left(
e^{-\pi i y\cdot B^{-1}A y}
 f\right)\right)_{B^{-1}}\right\|_{\fpq}\\
&\lesssim |\det
B|^{1/q-1/p-1/2}(\det(B^*B+I))^{1/2}|\det(I+iDB^{-1})|^{1/2}\\
&\qquad\qquad\quad\quad\quad\cdot\left\| \cF^{-1}\left(   e^{-\pi
i y\cdot B^{-1}A y}
 f\right)\right\|_{\fpq}\\
&\lesssim |\det B|^{1/q-1/p-1/2}(\det(B^*B+I))^{1/2}\det(I+iDB^{-1})|^{1/2}\\
&\qquad\qquad\quad\quad\quad\cdot\|    e^{-\pi i  y\cdot  B^{-1}A y} f\|_{\fqp}\\
&\lesssim \a(\mathcal{A},p,q)\| f\|_{\fqp}
\end{align*}
with $\a(\mathcal{A},p,q)$ given by \eqref{alfa}.

\noindent
 \emph{(ii)} In this case, we use formula \eqref{f5}.  Then, proceeding likewise the case \emph{(i)}, we majorize as
 follows:
\begin{align*}
 \|\mu(\A)f\|_{W(\cF
L^1,L^\infty)} &=|\det B|^{-1/2}\left\| e^{-\pi ix\cdot CA^{-1}x}
\left(e^{-\pi i y\cdot B^{-1}A y}\ast f\right)(A^{-1}
x)\right\|_{W(\cF
L^1,L^\infty)}\\
&\leq |\det B|^{-1/2} \|e^{-\pi i x\cdot CA^{-1}x}\|_{W(\cF
L^1,L^\infty)}\\
&\quad\quad\qquad\qquad\quad\,\,\,\,\,\,\,\,\,\,\,\,\cdot\left\|\left(
e^{-\pi i y\cdot B^{-1}A y}\ast f\right)_{A^{-1}}\right\|_{W(\cF
L^1,L^\infty)}\\
&\lesssim |\det B|^{-1/2}|\det
A|^{-1}(\det(A^*A+I))^{1/2}|\det(I+iCA^{-1})|^{1/2}\\
&\quad\quad\qquad\qquad\quad\,\,\,\,\,\,\,\,\,\,\,\,\cdot\left\|
 e^{-\pi i y\cdot B^{-1}A y}\ast
f\right\|_{W(\cF
L^1,L^\infty)}\\
&\lesssim \b(\mathcal{A})\| f\|_{W(\cF L^\infty,L^1)},
\end{align*}
where the last raw is due to
 \eqref{p2}, with $\b(\mathcal{A})$  defined
 in \eqref{beta}.
\end{proof}
\begin{proof}[Proof
of Theorem \ref{ttt1}] The proof uses the same arguments as in
Theorem \ref{t1}. Here, the estimate \eqref{dilAW0} is replaced by
\eqref{dilAW}. Besides,  the index relation \eqref{indexr} is
applied in the final step. In  details,
\begin{align*}
    \|\mu(\A)f\|_{\fpq}\!\! &\leq |\det B|^{-1/2}\|e^{-\pi i x\cdot DB^{-1}x}\|_{W(\cF
L^1,L^\infty)}\\
&\qquad\qquad\quad\quad\quad\cdot\left\|\left( \cF^{-1}\left(
e^{-\pi i y\cdot B^{-1}A y} f\right)\right)_{B^{-1}}\right\|_{\fpq}\\
&\lesssim \prod_{j=1}^d|\lambda_j|^{-1/2}|\det(I+iDB^{-1})(I+iB^{-1}A)|^{1/2}\\
&\quad\cdot\prod_{j=1}^d(\max\{1,|\lambda_j|^{-1}\})^{\mu_1(p',q')}
(\min\{1,|\lambda_j|^{-1}\})^{\mu_2(p',q')}\| f\|_{\fqp} \\
&= |\det(B+iD)(B+iA)|^{1/2}\\
&\quad\cdot \prod_{j=1}^d(\max\{1,|\lambda_j|\})^{\mu_1(p,q)-1/2}
(\min\{1,|\lambda_j|\})^{\mu_2(p,q)-1/2} \| f\|_{\fqp},
\end{align*}
that is case \emph{(i)}. Case
\emph{(ii)} indeed is not an
improvement of \eqref{f2} but is just
\eqref{f2} rephrased in terms of the
eigenvalues of $A$.
\end{proof}
\begin{remark}\rm The above theorems
 require the condition $\det B\not=0$. However, in some
special cases with $\det B=0$, the previous results can still be
used to obtain estimates between  Wiener amalgam spaces.
For example, if $\mathcal{A}=\begin{pmatrix} I&0\\
C&I\end{pmatrix}$, with $C=C^\ast$, then $\mu(\A)f(x)=\pm e^{-\pi i
C
 x\cdot x}f(x)$ (see \eqref{lower}), so that,
  for every $1\leq p,q\leq\infty$,
  Proposition \ref{p3} and the estimate
  \eqref{chirpnorm}
 give
 \[
 \|\mu(\A)f\|_{\fpq}\lesssim
 \prod_{j=1}^{d}(1+\lambda_j^2)^{1/4}\|f\|_{\fpq},
 \]
 where the $\lambda_j$'s are the
 eigenvalues of $C$
(incidentally, this estimate was already shown in
\cite{baoxiang,benyi,cordero}).
\end{remark}
\section{Applications to the Schr\"odinger
equation} In this section we apply the
previous results to the analysis of the
Cauchy problem of Schr\"odinger
equations with quadratic Hamiltonians,
i.e.
\begin{equation}\label{pc}
\begin{cases} i \displaystyle\frac{\partial
u}{\partial t} +H_\A u=0\\
u(0,x)=u_0(x),
\end{cases}
\end{equation}
where $H_\A$ is the Weyl quantization
of a quadratic form on the phase space
$\rdd$, defined from a matrix $\A$ in
the Lie algebra $\spdr$ of the
symplectic  group as follows (see
\cite{folland} and
\cite{degosson}).\par Any given matrix
$\A\in\spdr$ defines a quadratic form
$\mathcal{P}_\A(x,\xi)$ in $\rdd$ via
the formula
\[
\mathcal{P}_\A(x,\xi)=-\frac{1}{2}{}^t(x,\xi)\A\mathcal{J}(x,\xi),
\]
where, as usual,
$\mathcal{J}=\begin{pmatrix} 0&I\\
-I&0\end{pmatrix}$ (notice that
$\A\mathcal{J}$ is symmetric).
Explicitly, if
$\mathcal{A}=\begin{pmatrix}
A&B\\C&D\end{pmatrix}\in \spdr$ then
\begin{equation}\label{pa}
P_\A(x,\xi)=\frac{1}{2}\xi\cdot
B\xi-\xi\cdot A x-\frac{1}{2}x\cdot Cx.
\end{equation}
From the Weyl quantization, the
quadratic polynomial $P_\A$ in
\eqref{pa} corresponds to the Weyl
operator
 $\mathcal{P}_\A^w(D,X)$
 defined by
\[
2\pi
\mathcal{P}_\A^w(D,X)=-\frac{1}{4\pi}\sum_{j,k=1}^d
B_{j,k}\frac{\partial^2}{\partial
x_j\partial x_k}+i\sum_{j,k=1}^d
A_{j,k}
x_k\frac{\partial}{\partial{x_j}}+
\frac{i}{2}{\rm
Tr}(A)-\pi\sum_{j,k=1}^d C_{j,k}
x_j,x_k.
\]
The operator $H_\A:=2\pi
\mathcal{P}_\A^w(D,X)$ is called the
Hamiltonian operator.\par The evolution
operator for \eqref{pc} is related to
the metaplectic representation via the
following key formula
\[
e^{itH_\A}=\mu(e^{t\A}).
\]
Consequently,  Theorems \ref{t1} and
\ref{ttt1} can be used in the study of
fixed-time estimates for the solution
$u(t)=e^{itH_\A}u_0$ to \eqref{pc}.\par
As an  example, consider the matrix $\mathcal{A}=\begin{pmatrix} 0&B\\
0&0\end{pmatrix}\in\spdr$, with
$B=B^\ast$. Then the Hamiltonian
operator is
 $H_\A=-\frac{1}{4\pi}
 B\nabla\cdot\nabla$ and
 $e^{it\A}=\begin{pmatrix} I&
 tB\\ 0&I\end{pmatrix}\in Sp(d,\R)$.

Fix $t\not=0$. If $\det B\not=0$, and
$B$ has eigenvalues
  $\lambda_1,\dots,\lambda_d$,
 then the expression of $\beta'(e^{it\A})$ in \eqref{f2zb} is given by
 $$\beta'(e^{it\A})=2^{d/4}|\det tB|^{-1}|\det tB+iI)|^
 {1/2}=2^{d/4}\prod_{j=1}^d\left(
 \frac{1+t^2\lambda_j^2}{t^4\lambda_j^4}\right)^{1/4}.
 $$
Consequently, the fixed-time estimate
\eqref{f2z} is
 \[
 \|e^{it H_\A}
 f\|_{W(\Fur
 L^1,L^\infty)}\lesssim \prod_{j=1}^d \left(
 \frac{1+t^2\lambda_j^2}{t^4\lambda_j^4}\right)^{1/4}\|f
 \|_{W(\Fur
 L^\infty,L^1)},
 \]
 which generalizes the
 dispersive estimate in
 \cite{cordero}, corresponding
 to $B=I$.\par
 In the next two sections we
 present new fixed-time estimates, and
  also Strichartz estimates, in the
cases of the Hamiltonian
$H_\A=-\frac{1}{4\pi}\Delta+\pi|x|^2$
and
$H_\A=-\frac{1}{4\pi}\Delta-\pi|x|^2$.
\subsection{Schr\"odinger equation with
Hamiltonian
$H_\A=-\frac{1}{4\pi}\Delta+\pi|x|^2$}\
\par\noindent Here we
consider the Cauchy problem \eqref{pc}
with the Hamiltonian $H_\A$
corresponding to the matrix
$\A=\begin{pmatrix}
0&I\\-I&0\end{pmatrix}\in\spdr$, namely
$H_\A=-\frac{1}{4\pi}\Delta+\pi|x|^2$.
As a consequence of the estimates
proved in the previous section we
obtain the following fixed-time
estimates.
\begin{proposition}\label{p01}
For $2\leq r\leq\infty$, we have the
fixed-time estimates
\begin{equation}\label{d2}
\|e^{itH_\A} u_0\|_{\frpr}\lesssim|\sin
t|^{-2d\left(\frac{1}{2}-\frac{1}{r}\right)}
\|u_0\|_{\frrp}.
\end{equation}
\end{proposition}
\begin{proof}
The symplectic matrix $e^{t\A}$ reveals
to be $e^{t\A}=\begin{pmatrix}(\cos
t)I&(\sin t)I\\(-\sin t)I&(\cos
t)I\end{pmatrix}.$

First, using the estimate \eqref{f2z}
we get
\begin{equation}\label{d1}
\|e^{itH_\A} u_0\|_{\fui}\lesssim|\sin
t|^{-d} |\cos
t|^{-\frac{3}{2}d}\|u_0\|_{\fiu}.
\end{equation}
On the other hand, the estimate
\eqref{f1z}, for $p=1, q=\infty$, reads
\begin{equation}\label{d1w}
\|e^{itH_\A} u_0\|_{\fui}\lesssim|\sin
t|^{-5d/2}\|u_0\|_{\fiu}.
\end{equation}
 Since $\min\{|\sin
t|^{-d} |\cos t|^{-\frac{3}{2}d},|\sin
t|^{-5d/2}\} \asymp|\sin t|^{-d}$, we
obtain \eqref{d2} for $r=\infty$, which
is the dispersive estimate.\par The
estimates \eqref{d2} for $2\leq r\leq
\infty$ follow by complex interpolation
from the dispersive estimate and the
$L^2-L^2$ estimate
\begin{equation}\label{r01}
\|e^{itH_\A}f\|_{L^2}= \|f\|_{L^2}.
\end{equation}

\end{proof}\par\noindent

The Strichartz estimates for the
solutions to \eqref{pc} are detailed as
follows.

\begin{theorem}\label{prima} Let
$T>0$ and $4<q,\tilde{q}\leq\infty$,
$2\leq r,\tilde{r}\leq\infty$, such
that
\begin{equation}\label{0000}\frac{2}{q}+\frac{d}{r}=\frac{d}{2},
\end{equation} and
similarly for $\tilde{q},\tilde{r}$.
Then we have the homogeneous Strichartz
estimates
\begin{equation}\label{hom}\|e^{itH_{\A}}
u_0\|_{L^{{q}/{2}}([0,T]) W(\Fur
L^{r^\prime},L^r)_x}\lesssim
\|u_0\|_{L^2_x},
\end{equation}
the dual homogeneous Strichartz
estimates
\begin{equation}\label{dh}
\|\int_0^T e^{-isH_{\A}}
F(s)\,ds\|_{L^2}\lesssim
\|F\|_{L^{(\tilde{q}/{2})^\prime}([0,T])
W(\Fur
L^{\tilde{r}},L^{\tilde{r}^\prime})_x},
\end{equation}
and the retarded Strichartz estimates
\begin{equation}\label{ret}
\|\int_{0\leq s<t} e^{i(t-s)H_{\A}}
F(s)\,ds\|_{L^{q/2}([0,T])W(\Fur
L^{r^\prime},L^r)_x}
\lesssim\|F\|_{L^{(\tilde{q}/{2})^\prime}([0,T])
W(\Fur
L^{\tilde{r}},L^{\tilde{r}^\prime})_x}.
\end{equation}
Consider then the endpoint
$P:=(4,2d/(d-1))$. For $(q,r)=P$,
$d>1$, we have
\begin{equation}\label{hombis}
\|e^{itH_{\A}} u_0\|_{L^{2}([0,T])
W(\Fur L^{r^\prime,2},L^r)_x}\lesssim
\|u_0\|_{L^2_x},
\end{equation}
\begin{equation}\label{dhbis}
\|\int_0^T e^{-isH_{\A}}
F(s)\,ds\|_{L^2}\lesssim
\|F\|_{L^{2}([0,T]) W(\Fur
L^{r,2},L^{r^\prime})_x}.
\end{equation}
The retarded estimates \eqref{ret}
still hold with $(q,r)$ satisfying
\eqref{0000}, $q>4,r\geq2$,
$(\tilde{q},\tilde{r})=P$, if one
replaces $\Fur L^{\tilde{r}'}$ by $\Fur
L^{\tilde{r}',2}$. Similarly it holds
for $(q,r)=P$ and
$(\tilde{q},\tilde{r})\not=P$ as above
if one replaces $\Fur L^{r'}$ by $\Fur
L^{r',2}$. It holds for both
$(p,r)=(\tilde{p},\tilde{r})=P$ if one
replaces $\Fur L^{r'}$ by $\Fur
L^{r',2}$ and $\Fur L^{\tilde{r}'}$ by
$\Fur L^{\tilde{r}',2}$.
\end{theorem}
In the previous theorem the bounds may
depend on $T$.
\begin{proof}
The arguments are essentially the ones
in \cite{cordero, keel}. For the
convenience of the reader, we present
 the guidelines of the proof.\par Due to the property
group of the evolution operator
$e^{itH_\A}$, we can limit ourselves to
the case $T=1$. Indeed, observe that,
if \eqref{hom} holds for
 a given $T>0$, it holds for any $0< T'\leq T$
  as well, so that it
  suffices to
  prove \eqref{hom}
 for $T=N$ integer. Since
\[
\|e^{itH_\A}u_0\|^{\frac{q}{2}}_{L^{q/2}
([0,N]) W(\Fur
L^{r^\prime},L^r)_x}=\sum_{k=0}^{N-1}\|
e^{itH_\A}e^{ikH_\A}u_0\|^{\frac{q}{2}}_{
L^{q/2}([0,1])W(\Fur
L^{r^\prime},L^r)_x},
\]
the  $T=N$ case is reduced to the $T=1$
case  by using \eqref{hom} for $T=1$
and the conservation law \eqref{r01}.
 The other estimates
can be treated analogously. Whence from
now on $T=1$.\par
 Consider
first the non-endpoint case. Set
$U(t)=\chi_{[0,1]}(t)e^{itH_\A}$. For
$2\leq r\leq\infty$, using relation
\eqref{d2}, we get
\begin{equation}\label{d3}
\|U(t)(U(s))^\ast
f\|_{\frpr}\lesssim|t-s|^{-2d\left(\frac{1}{2}-\frac{1}{r}\right)}\|f\|_{\frrp}.
\end{equation}
By the $TT^\ast$ method\footnote{This
duality argument is generally
established for $L^p$ spaces. Its use
 for Wiener amalgam spaces is similarly
justified thanks to the duality defined
 by the H\"older-type
inequality \cite{cordero}:
\[
|\langle F,G\rangle_{L^2_t L^2_x}|\leq
\|F\|_{W(L^s,L^q)_t W(\Fur
L^{r^\prime},L^r))_x}\|G\|_{W(L^{s^\prime},L^{q^\prime})_t
W(\Fur L^{r},L^{r^\prime})_x}.\]} (see,
e.g., \cite[Lemma 2.1]{GinibreVelo92}
or \cite[page 353]{stein}) the estimate
\eqref{hom} is equivalent to
\begin{equation}\label{d4}
\|\int U(t) (U(s))^\ast
F(s)\,ds\|_{L^{q/2}_tW(\Fur
L^{r^\prime},L^r)_x}
\lesssim\|F\|_{L^{({q}/{2})^\prime}_t
W(\Fur
L^{\tilde{r}},L^{{r}^\prime})_x}.
\end{equation}
The estimate above is attained by
applying   Minkowski's inequality and
the Hardy-Littlewood-Sobolev inequality
\eqref{conv1} to the estimate
\eqref{d3}. The dual homogeneous
estimates \eqref{dh} follow by duality.
Finally, the retarded estimates
\eqref{ret}, with $(1/q,1/r)$,
$(1/\tilde{q}, 1/\tilde{r})$ and
$(1/\infty,1/2)$ collinear, follow by
complex interpolation from the three
cases $(\tilde{q},\tilde{r})=(q,r)$,
$(q,r)=(\infty,2)$ and
$(\tilde{q},\tilde{r})=(\infty,2)$,
which in turns are a consequence,  of
\eqref{d4} (with $\chi_{s<t}F$ in place
of $F$), \eqref{dh} (with $\chi_{s<t}F$
in place of $F$) and the duality
argument, respectively.\par We are left
to the endpoint case:
 $(q,r)=(2,2d/(d-1))$. The estimate \eqref{hombis} is
equivalent to the bilinear estimate
\[
|\iint \langle (U(s))^\ast F(s),
(U(t))^\ast
G(t)\rangle\,ds\,dt|\lesssim
\|F\|_{L^2_t W(\Fur
L^{{r,2}},L^{{r}^\prime})_x}
\|G\|_{L^2_t W(\Fur
L^{{r,2}},L^{{r}^\prime})_x}.
\]
By symmetry, it is enough to prove
\begin{equation}\label{aim}
|T(F,G)|\lesssim \|F\|_{L^2_t W(\Fur
L^{{r},2},L^{{r}^\prime})_x}
\|G\|_{L^2_t W(\Fur
L^{{r},2},L^{{r}^\prime})_x},
\end{equation}
where
\[
T(F,G)=\iint_{s<t} \langle (U(s))^\ast
F(s), (U(t))^\ast G(t)\rangle\,ds\,dt.
\]
To this aim,  $T(F,G)$ is  decomposed
dyadically
 as
$T=\sum_{j\in \mathbb{Z}} T_j$, with
\begin{equation}\label{tj}
T_j(F,G)=\iint_{t-2^{j+1}<s\leq t-2^j}
(U(s))^\ast F(s), (U(t))^\ast
G(t)\rangle\,ds\,dt.
\end{equation}
By resorting on \eqref{dh} one can
prove exactly as in \cite[Lemma
4.1]{keel} the following estimates:
\begin{equation}\label{lammasecondo}
|T_j(F,G)|\lesssim
2^{-j\beta(a,b)}\|F\|_{L^2_t W(\Fur
L^a,L^{a^\prime})}\|G\|_{L^2_t W(\Fur
L^b,L^{b^\prime})},
\end{equation}
for $(1/a,1/b)$ in a neighborhood of
$(1/r,1/r)$, with
$\beta(a,b)=d-1-\frac{d}{a}-\frac{d}{b}$.\par
The estimate \eqref{aim} is achieved by
means of a  real interpolation result,
detailed in \cite[Lemma 6.1]{keel}, and
applied to the vector-valued bilinear
operator $T=(T_j)_{j\in\mathbb{Z}}$.
Here, however, we must observe that, if
$A_k=L^2_tW(\Fur
L^{a_k},L^{{a_k}^\prime})_x$, $k=0,1$,
and $\theta_0$ fulfills
$1/r=(1-\theta_0)/a_0+\theta_0/a_1,$
then
 \[
L^2_tW(\Fur L^{r,2},L^{{r}^\prime})_x
\subset (A_0,A_1)_{\theta_0,2}. \]
 The above inclusion follows by
\cite[Theorem 1.18.4, page
129]{triebel} (with $p=p_0=p_1=2$) and
Proposition \ref{inter9}. This gives
\eqref{hombis} and \eqref{dhbis}.\par
Consider now the endpoint retarded
estimates. The case
$(\tilde{q},\tilde{r})=(q,r)=P$ is
exactly \eqref{aim}. The case
$(\tilde{q},\tilde{r})=P$,
$(q,r)\not=P$,  can be obtained by a
repeated use of H\"older's inequality
to interpolate from the case
$(\tilde{q},\tilde{r})=(q,r)=P$ and the
case $(\tilde{q},\tilde{r})=P$,
$(q,r)=(\infty,2)$ (that is clear from
\eqref{dhbis}). Finally,  the retarded
estimate in the case $(q,r)=P$,
$(\tilde{q},\tilde{r})\not=P$, follows
by applying the arguments above to the
adjoint operator $ G\mapsto \int _{t>s}
(U(t))^\ast U(s) G(t)\,dt$, which gives
the dual estimate.
\end{proof}
\subsection{Schr\"odinger equation with
Hamiltonian
$H_\A=-\frac{1}{4\pi}\Delta-\pi|x|^2$}\
\par\noindent
The   Hamiltonian operator
$H_\A=-\frac{1}{4\pi}\Delta-\pi|x|^2$
corresponds to the matrix
$\A=\begin{pmatrix}
0&I\\I&0\end{pmatrix}\in\spdr$. In this
case, $e^{t\A}=\begin{pmatrix}(\cosh
t)I&(\sinh t)I\\(\sinh t)I&(\cosh
t)I\end{pmatrix}\in Sp(d,\R)$.\par
Fixed-time estimates for $H_\A$ are as
follows.
\begin{proposition} For
$2\leq r\leq\infty$,
\begin{equation}\label{d20}
\|e^{itH_\A}
u_0\|_{\frpr}\lesssim\left(\frac{1+|\sinh
t|}{\sinh^2
t}\right)^{d\left(\frac{1}{2}-\frac{1}{r}\right)}
\|u_0\|_{\frrp}.
\end{equation}
\end{proposition}
\begin{proof}
The estimate \eqref{f2z} yields the
dispersive estimate
\begin{equation}\label{d10}
\|e^{itH_\A}
u_0\|_{\fui}\lesssim\left(\frac{1+|\sinh
t|}{\sinh^2 t}\right)^{\frac{d}{2}}
\|u_0\|_{\fiu}.
\end{equation}
(Observe that \eqref{f1z}, with $p=1$,
$q=\infty$, gives a bound worse than
\eqref{d10}).\par The estimates
\eqref{d20} follow by complex
interpolation between the dispersive
estimate \eqref{d10} and the
conservation law \eqref{r01}.
\end{proof}\par\noindent
We can now establish the corresponding
Strichartz estimates.
\begin{theorem}\label{seconda}
Let $4<q,\tilde{q}\leq\infty$, $2\leq
r,\tilde{r}\leq\infty$, such that
\begin{equation}\label{0001}\frac{2}{q}+\frac{d}{r}=\frac{d}{2},
\end{equation} and
similarly for $\tilde{q},\tilde{r}$.
Then we have the homogeneous Strichartz
estimates
\begin{equation}\label{hom0}\|e^{itH_\A}
u_0\|_{W(L^{{q}/{2}},L^2)_t W(\Fur
L^{r^\prime},L^r)_x}\lesssim
\|u_0\|_{L^2_x},
\end{equation}
the dual homogeneous Strichartz
estimates
\begin{equation}\label{dh0}
\|\int e^{-isH_\A}
F(s)\,ds\|_{L^2}\lesssim
\|F\|_{W(L^{(\tilde{q}/{2})^\prime},L^{
2})_t W(\Fur
L^{\tilde{r}},L^{\tilde{r}^\prime})_x},
\end{equation}
and the retarded Strichartz estimates
\begin{equation}\label{ret0}
\|\int_{s<t} e^{i(t-s)H_\A}
F(s)\,ds\|_{W(L^{q/2},L^2)_t W(\Fur
L^{r^\prime},L^r)_x}
\lesssim\|F\|_{W(L^{(\tilde{q}/{2})^\prime},
L^{2})_t W(\Fur
L^{\tilde{r}},L^{\tilde{r}^\prime})_x}.
\end{equation}
Consider then the endpoint
$P:=(4,2d/(d-1))$. For $(q,r)=P$,
$d>1$, we have
\begin{equation}\label{hombis0}
\|e^{itH_\A} u_0\|_{L^2_t W(\Fur
L^{r^\prime,2},L^r)_x}\lesssim
\|u_0\|_{L^2_x},
\end{equation}
\begin{equation}\label{dhbis0}
\|\int e^{-isH_\A}
F(s)\,ds\|_{L^2}\lesssim \|F\|_{L^{2}_t
W(\Fur L^{r,2},L^{r^\prime})_x}.
\end{equation}
The retarded estimates \eqref{ret0}
still hold with $(q,r)$ satisfying
\eqref{0001}, $q>4,r\geq2$,
$(\tilde{q},\tilde{r})=P$, if one
replaces $\Fur L^{\tilde{r}'}$ by $\Fur
L^{\tilde{r}',2}$. Similarly it holds
for $(q,r)=P$ and
$(\tilde{q},\tilde{r})\not=P$ as above
if one replaces $\Fur L^{r'}$ by $\Fur
L^{r',2}$. It holds for both
$(p,r)=(\tilde{p},\tilde{r})=P$ if one
replaces $\Fur L^{r'}$ by $\Fur
L^{r',2}$ and $\Fur L^{\tilde{r}'}$ by
$\Fur L^{\tilde{r}',2}$.
\end{theorem}
\begin{proof}
Let us first prove \eqref{hom0}. By the
$TT^\ast$ method it suffices to prove
\begin{equation}\label{w04}
\|\int
e^{i(t-s)H_\A}F(s)\,ds\|_{W({L^{q/2}},{L^2})_t
W(\Fur L^{r^\prime},L^r)_x}\lesssim
\|F\|_{W({L^{\left(q/2\right)^\prime}},
L^{2})_t W(\Fur L^r,L^{r^\prime})_x}.
\end{equation}
For $0<\alpha<1/2$, let
$\phi_\alpha(t)=|\sinh
t|^{-\alpha}+|\sinh t|^{-2\alpha}$,
$t\in\mathbb{R}$, $t\not=0$. A direct
computation shows that $\phi_\alpha\in
W(L^{1/(2\alpha),\infty},L^1)$. Since
 $L^1\ast L^2\hookrightarrow L^2$
(Young's Inequality) and
$L^{\left(\frac{1}{\alpha}\right)'}\ast
L^{\frac{1}{2\alpha},\infty}\hookrightarrow
L^{\frac{1}{\alpha}}$ (Proposition
\ref{convlor}), Lemma \ref{WA} $(i)$
gives the convolution relation
\begin{equation}\label{w03}
\|F\ast
\phi_\alpha\|_{W(L^{1/\alpha},L^{2/\alpha})}\lesssim
\|F\|_{W(L^{{(1/\alpha)}^\prime},L^{(2/\alpha)^\prime})}.
\end{equation}
Fix now $\alpha=d(1/2-1/r)=2/q$; then,
by \eqref{d20},  \eqref{w03} and
Minkowski's Inequality,
\begin{align*}
\|\int e^{i(t-s)H_\A}F(s)\,ds &
\|_{W(L^{q/2},{L^2})_t W(\Fur
L^{r^\prime},L^r)_x}\\
&\leq \left\|\int\|e^{i(t-s)H_\A}
F(s)\|_{W(\Fur
L^{r^\prime},L^r)_x}\,ds\right\|_{W({L^{q/2}},
{L^2})_t}\\
&\lesssim \|\|F(t)\|_{W(\Fur
L^{r^\prime},L^r)_x}\ast
\phi_\alpha(t)\|_{W(L^{q/2},{L^2})_t}\\
&\lesssim \|F\|_{W(L^{(q/2)^\prime},
L^{2})_t W(\Fur L^r,L^{r^\prime})_x}.
\end{align*}
This proves \eqref{w04} and whence
\eqref{hom0}.  The estimate \eqref{dh0}
follows from \eqref{hom0} by duality.
The proof of \eqref{ret0} is analogous
to \eqref{ret} in Theorem
\ref{prima}.\par For the endpoint case
one can repeat essentially verbatim the
arguments in the proof of Theorem
\ref{prima}, upon setting
$U(t)=e^{itH_\A}$. To avoid
repetitions, we omit the details (see
also the proof of \cite[Theorem
1.2]{cordero}).
\end{proof}
\subsection{Comparison with
the classical estimates in Lebesgue
spaces} Here we compare the above
estimates with the classical ones
between Lebesgue spaces. For the
convenience of the reader we recall the
following very general result by Keel
and Tao \cite[Theorem 1.2]{keel}.\par
Given $\sigma>0$, we say that an
exponent pair $(q,r)$ is {\it sharp
$\sigma$-admissible} if
$1/q+\sigma/r=\sigma/2$,
$q\geq2,r\geq2$,
$(q,r,\sigma)\not=(2,\infty,1)$.
\begin{theorem}\label{tkt}
Let $(X,\mathcal{S},\mu)$ be a
$\sigma$-finite measured space, and
$U:\R\to B(L^2(X,\mathcal{S},\mu))$ be
a weakly measurable map satisfying, for
some $\sigma>0$,
\[
\|U(t)f\|_{L^2}\lesssim
\|u\|_{L^2},\quad t\in\R,
 \]
and
\[
\|U(s)U(t)^\ast f\|_{L^\infty}\lesssim
|t-s|^{-\sigma}\|f\|_{L^1},\quad
t,s\in\R.
\]
Then for every sharp
$\sigma$-admissible pairs $(q,r)$,
$(\tilde{q},\tilde{r})$, one has
\[
\|U(t) f\|_{L^q_t L^r_x}\lesssim
\|f\|_{L^2},
\]
\[
\|\int U(s)^\ast
F(s)\,ds\|_{L^2}\lesssim
\|F\|_{L^{q'}_t L^{r'}_x},
\]
\[
\|\int _{s<t} U(t) U(s)^\ast
F(s)\,ds\|_{L^q_t L^r_x}\lesssim
\|F\|_{L^{\tilde{q}'}L^{\tilde{r}'}}.
\]
\end{theorem}
First we fix the attention to the case
of the Hamiltonian
$H_\A=-\frac{1}{4\pi}\Delta+\pi|x|^2$.
One has the following explicit formula
for $e^{it H_\A}u_0= \mu(e^{it\A})u_0$
in \eqref{f3}:
\[
e^{itH_\A}u_0=i^{d/2}(\sin
t)^{-d/2}\int e^{-\pi i ({\rm cotg} t)
(|x|^2+|y|^2)+2\pi i ({\rm cosec} t)
y\cdot x}u_0(y)\,dy.
\]
Hereby it follows immediately the
dispersive estimate
\begin{equation}\label{ll1}
\|e^{itH_\A} u_0\|_{L^\infty}\leq |\sin
t|^{-d/2}\|u\|_{L^1}.
\end{equation}
Notice that \eqref{ddd2} (i.e.
\eqref{d2} with $r=\infty$) represents
an improvement of \eqref{ll1} for every
fixed $t\not=0$, since
$L^1\hookrightarrow W(\Fur
L^\infty,L^1)$ and $W(\Fur
L^1,L^\infty)\hookrightarrow L^\infty$.
However, as might be expected, the
bound on the norm in \eqref{ddd2}
becomes worse than that in \eqref{ll1}
as $t\to k\pi$, $k\in\mathbb{Z}$.\par
As a consequence of \eqref{ll1},
Theorem \ref{tkt} with
$U(t)=e^{itH_\A}\chi_{[0,1]}(t)$ and
$\sigma=d/2$, and the group property of
the operator $e^{itH_\A}$ (as in the
proof of Theorem \ref{prima} above) one
deduce, for example, the homogeneous
Strichartz estimate
\begin{equation}\label{lmaa}
\|e^{itH_{\A}} u_0\|_{L^{{q}}([0,T])
L^r_x}\lesssim \|u_0\|_{L^2_x},
\end{equation}
for every pair $(q,r)$ satisfying
$2/q+d/r=d/2$, $q\geq 2,r\geq2$,
$(q,r,d)\not=(2,\infty,2)$. These
estimates were also
 obtained recently in
\cite{koch} by different methods.\par
Hence, one sees that \eqref{hom}
predicts, for the solution to
\eqref{pc}, a better local spatial
regularity than \eqref{lmaa}, but just
after averaging on $[0,T]$ by the
$L^{q/2}$ norm, which is smaller than
the $L^q$ norm.\par\medskip We now
consider the case of the Hamiltonian
$H_\A=-\frac{1}{4\pi}\Delta-\pi|x|^2$.\\
The dispersive estimate here reads
\begin{equation}\label{pol}
\|e^{itH_\A}u_0\|_{L^{\infty}(\R^d)}\leq
|\sinh t |^{-d/2}
\|u_0\|_{L^{1}(\bR^d)}.
\end{equation}
This estimate follows immediately from
the
 explicit
 expression of
$e^{it H_\A}u_0= \mu(e^{it\A})u_0$ in
\eqref{f3}:
\begin{equation*} e^{it H_\A}u_0=i^{d/2}
(\sinh t)^{-d/2}\int e^{-\pi i ({\rm
cotgh} t) (|x|^2+|y|^2)+2\pi i ({\rm
cosech} t) y\cdot x} u_0(y)\,dy.
\end{equation*}
The corresponding  Strichartz estimates
 between the Lebesgue spaces
read
\begin{equation}\label{S1}
\|e^{it H_\A} u_0\|_{L^q_t
L^r_x}\lesssim \|u_0\|_{L^2_x},
\end{equation}
for $q\geq2$, $r\geq 2$, with
$2/q+d/r=d/2$,
$(q,r,d)\not=(2,\infty,2)$. These
estimates are the issues of Theorem
\ref{tkt} with $U(t)=e^{itH_\A}$, and
the dispersive estimate \eqref{pol}
(indeed, $|\sinh
t|^{-d/2}\leq|t|^{-d/2}$). These
estimates are to be compared with
\eqref{d20} (with $r=\infty$) and
\eqref{hom0} respectively.\par
 One can
do the same remarks as in the previous
case. In addition here one should
observe that \eqref{hom0} displays a
better time decay at infinity than the
classical one ($L^2$ instead of $L^r$),
for a norm, $\|u(t,\cdot)\|_{W(\Fur
L^{r'},L^r)}$, which is even bigger
than $L^r$. Notice however that our
range of exponents is restricted to
$q\geq 4$.

\section*{Acknowledgements}
The authors would like to thank Prof. Maurice de Gosson and Luigi
Rodino for fruitful conversations and comments.

\vskip0.5truecm

\end{document}